\documentclass[11pt]{article}
\usepackage{amsfonts}

\usepackage{amsfonts, amsmath, amssymb}
\usepackage{amssymb,amsfonts,amsmath,latexsym, epsfig,cite, psfrag,eepic,colordvi}
\usepackage{amscd,graphics}
\usepackage{lineno}

\newcommand{\qed}{\hfill $\Box $}
\newcommand{\pf}{\noindent {\bf Proof.} }

\newtheorem{theorem}{Theorem}[section]
\newtheorem{lemma}[theorem]{Lemma}
\newtheorem{coro}[theorem]{Corollary}

\newtheorem{conjecture}[theorem]{Conjecture}

\begin{document}

\title{Vertex-Coloring 2-Edge-Weighting of Graphs
\thanks{This work is supported in part by 973 Project of Ministry of Science and
Technology of China and the Natural Sciences and Engineering
Research Council of Canada.}}
\author{Hongliang Lu\,\textsuperscript{a}, \
Qinglin Yu\,\textsuperscript{b} and Cun-Quan
Zhang\,\textsuperscript{c}
\\ {\small \textsuperscript{a}Department of Mathematics}
\\ {\small  Xi'an Jiaotong University, Xi'an,  China
}
\\ {\small \textsuperscript{b}Department of Mathematics and Statistics}
\\ {\small Thompson Rivers University, Kamloops, BC, Canada}
\\ {\small \textsuperscript{c}Department of Mathematics}
\\ {\small West Virginia University, Morgantown, WV, USA}
}

\date{}

\maketitle

\begin{abstract}
A $k$-{\it edge-weighting} $w$ of a graph $G$ is an assignment of an
integer weight, $w(e)\in \{1,\dots , k\}$, to each edge $e$. An edge
weighting naturally induces a vertex coloring $c$ by defining
$c(u)=\sum_{u\sim e} w(e)$ for every $u \in V(G)$. A
$k$-edge-weighting of a graph $G$ is \emph{vertex-coloring} if the
induced coloring $c$ is proper,  i.e., $c(u) \neq c(v)$  for any edge
$uv \in E(G)$.

Given a graph $G$ and a vertex coloring $c_0$, does there exist an
edge-weighting such that the induced vertex coloring is $c_0$? We
investigate this problem by considering edge-weightings defined on
an abelian group.

It was proved that every 3-colorable graph admits a vertex-coloring
$3$-edge-weighting \cite{KLT}. Does every 2-colorable graph (i.e.,
bipartite graphs) admit a vertex-coloring 2-edge-weighting?  We
obtain several simple sufficient conditions for graphs to be
vertex-coloring 2-edge-weighting.  In particular, we show that
3-connected bipartite graphs admit vertex-coloring 2-edge-weighting.
\\
[2mm] \textbf{Keywords}: edge-weighting; vertex-coloring; $3$-connected bipartite graph.\\
[2mm] {\bf AMS subject classification (2000)}: 05C15.
\end{abstract}

\section{Introduction}

In this paper, we consider only finite, undirected and simple
connected graphs. For a vertex $v$ of a graph $G=(V,E)$, $N_G(v)$
denotes the set of vertices which are adjacent to $v$. If $v\in
V(G)$ and $e\in E(G)$, we use $v\sim e$ to denote that $v$ is an
end-vertex of $e$, $\omega(G)$ denotes the number of connected
components of $G$. An $k$-{\em vertex coloring} $c$ of $G$ is an
assignment of $k$ integers, $1,2,\dots, k$, to the vertices of $G$,
the color of a vertex $v$ is denoted by $c(v)$. The coloring is {\em
proper} if no two distinct adjacent vertices share the same color. A
graph $G$ is $k$-{\em colorable} if $G$ has a proper $k$-vertex
coloring. The {\em chromatic number} $\chi(G)$ is the minimum number
$r$ such that $G$ is $r$-colorable. Notation and terminology that is
 not defined here may be found in \cite{BB}.

A $k$-{\it edge-weighting} $w$ of a graph $G$ is an assignment of an
integer weight $w(e)\in \{1,\dots , k\}$ to each edge $e$ of $G$. An
edge weighting naturally induces a vertex coloring $c(u)$ by
defining $c(u)=\sum_{u\sim e} w(e)$ for every $u \in V(G)$. An
$k$-edge-weighting of a graph $G$ is \textit{vertex-coloring} if for
every edge $e=uv$, $c(u)\neq c(v)$ and then we say $G$ admitting a
\emph{vertex-coloring $k$-edge-weighting}. Moreover, we introduce a
concept, which  is different from the concept discussed here but
similar enough.  A multigraph is \emph{irregular} if no two vertex
degrees are equal. A multigraph can be viewed as a weighted graph
with nonnegative-integer weights on the edges. The degree of a
vertex in a weighted graph is the sum of the incident weights.
Chartrand et al. \cite{CJLORS} defined the \emph{irregularity
strength} of a graph $G$, written $s(G)$, to be the minimum of the
maximum edge weight in an irregular multigraph with underlying graph
$G$.

If a graph has an edge as a component, clearly it can not have a
vertex-coloring $k$-edge-weighting. So in this paper, we only
consider graphs without $K_2$ component and refer such graphs as
\textit{nice graphs}.

In \cite{KLT}, Karo\'{n}ski, {\L}uczak and Thomason initiated the
study of vertex-coloring  $k$-edge-weighting and  they brought
forward a conjecture as following.

\begin{conjecture}{\rm{(1-2-3-Conjecture)}}
Every nice graph admits a vertex-coloring 3-edge-weighting.
\end{conjecture}

Furthermore, they proved that the conjecture holds for 3-colorable
graphs (see Theorem 1 in \cite{KLT}). For other graphs,
Addario-Berry \textit{et al.} \cite{ADMRT} showed that every nice
graph admits a vertex-coloring 30-edge-weighting. Addario-Berry,
Dalal and Reed \cite{ADR} improved the number of integers required
to 16.  Later, Wang and Yu \cite{WY} improved this bound to 13.
Recently, Kalkowski, Karo\'{n}ski and Pfender \cite{KKP} showed that
every nice graph admits a vertex-coloring 5-edge-weighting, which is
a great leap towards the 1-2-3-Conjecture.

In this paper, we focus on vertex-coloring 2-edge-weighting. In
Section 2, we present several new results about vertex-coloring
2-edge-weighting.

Besides the existence problem of vertex-coloring $k$-edge-weighting,
a natural question to ask is that given a graph $G$ and a vertex
coloring $c_0$, can we realize the coloring $c_0$ by a
$k$-edge-weighting, i.e., does there exist  an edge-weighting such
that the induced vertex coloring is $c_0$? For general graphs, it is
not easy to find such an edge-weighting. However, if restricting
edge weights to an abelian group, we obtain a neat positive answer
for this even for a non-proper coloring $c_0$. In Section 3, we show
that every 3-connected nice bipartite graph admits a vertex-coloring
2-edge-weighting.

\section{Vertex-coloring $2$-edge-weighting}

For a graph $G$, there is a close relationship between
2-edge-weightings and graph factors. Namely, a 2-edge-weighting
problem is equivalent to finding a special factor of graphs (see
\cite{ADMRT,ADR}). So to find spanning subgraphs with pre-specified
degree is an important part of edge-weighting. We shall use some of
these results in our proofs.

\begin{lemma}{\rm{(Addario-Berry, Dalal and Reed, \cite{ADR})}}\label{Berry}
Given a graph $G=(V,E)$, if for all $v\in V$, there are integers
$a_{v}^{-}$, $a_{v}^{+}$ such that $a_{v}^{-}\leq \lfloor
\frac{1}{2}d(v)\rfloor\leq a_{v}^{+}<d(v)$, and
\begin{align*}
a_{v}^{+}\leq \min\{\frac{1}{2}(d(v)+a_{v}^{-})+1,
2(a_{v}^{-}+1)+1\},
\end{align*}
then there exists a spanning subgraph $H$ of $G$ such that
$d_{H}(v)\in \{a_{v}^{-},a_{v}^{-}+1,a_{v}^{+},a_{v}^{+}+1\}$.
\end{lemma}


Given an arbitrary vertex coloring $c_0$, we want to find an
edge-weighting such that the induced vertex coloring is $c_0$? Under
a weak condition, the next two theorems show that there exists an
edge-weighting from an abelian group to $E(G)$ to induce $c_0$ for
bipartite and non-bipartite graphs respectively.

\begin{theorem}\label{group2}
Let $G$ be a   non-bipartite graph and $\Gamma=\{g_1,g_2,\dots,
g_k\}$ be a finite abelian group, where $k=|\Gamma|$. Let $c_0$ be
any $k$-vertex coloring of $G$ with color classes $\{U_{1},\dots,
U_k\}$, where $|U_i|=n_i$ $(1\leq i\leq k)$. If there exists an
element $h\in \Gamma$ such that $n_{1}g_{1}+\dots+n_{k}g_{k}=2h$,
then there is an edge-weighting with the elements of $\Gamma$ such
that the induced vertex coloring is $c_0$.
\end{theorem}

\pf Let $c_0$ be any  $k$-vertex coloring with vertex partition
$\{U_{1},\dots, U_k \}$, where every element in $U_{i}$ is colored
with $g_{i}$ $(1\leq i\leq k)$ such that
$n_{1}g_{1}+\dots+n_{k}g_{k}=2h$.

Assign one edge with weight $h$ and the rest with zero, so the sum
of vertex colors is $2h$. We now adjust this initial weighting,
while maintaining the sum of vertex weights, until all the vertices
in $U_{i}$  have color $g_i$ $(1\leq i\leq k)$. Suppose there exists
a vertex $u \in U_i$ with the wrong color $g \neq g_i$. Since
$n_{1}g_{1}+\dots+n_{k}g_{k}=2h$, there must be another vertex $v\in
V(G)$ whose color is also wrong. Since $G$ is non-bipartite, we can
choose a walk of even length from $u$ to $v$, which is always
possible since $k \geq 3$. Traverse this walk, adding $g_i-g, g-g_i,
g_i-g,\dots$ alternately to the edges as they are encountered. This
operation maintains the sum of vertex weights, leaves the colors of
all but $u$ and $v$ unchanged, and yields one more vertex of correct
color. Hence, repeated applications give the desired weighting. \qed

\begin{theorem}\label{group3}
Let $G$ be a nice bipartite graph and $Z_{2}=\{0,1\}$. Let $c_0$ be
any $2$-vertex coloring of $G$ with color classes $\{U_{0}, U_1\}$,
where $|U_i|=n_i$ $(0\leq i\leq 1)$ such that $c_{0}(U_i)=i$,  for
$i=0,1$. If $n_{1}$ is even, then there exists an edge-weighting
with the elements of $Z_{2}$ such that the induced vertex coloring
is $c_0$.
\end{theorem}

\pf  Let $g_1=0$ and $g_2=1$. If there is a vertex $u$ of color
$g_i$ with the wrong color $g\neq g_i$ and since $n_{2}$ is even,
then there must be another vertex $v\in V(G)$ whose color is also
wrong. Since $G$ is connected, then there is a path from $u$ to $v$.
Traverse this walk, adding $1,1,1,\dots$ to the edges as they are
encountered. This operation always maintains the sum of vertex
colors, leaves the colors of all but $u$ and $v$ unchanged, and
yields one more vertex of correct weight. \qed

%
%
%

Note that in Theorem \ref{group2}, the given vertex-coloring $c_0$
can be either a proper or an improper coloring.
\vspace{3mm}

\textbf{Remark:} The edge-weighting problem on groups has been
studied by Karo\'{n}ski, {\L}uczak and Thomason in \cite{KLT}. They
proved that for each $|\Gamma|$-colorable graph $G$,  there exists
an edge-weighting with the elements of $\Gamma$ such that the
induced vertex-coloring is proper. Our proof of Theorems
\ref{group2} and \ref{group3} are modifications of the result.

It was proved in \cite{KLT} that every $3$-colorable graph has a
vertex-coloring $3$-edge-weighting. A natural question to ask is
that whether every $2$-colorable graph has a vertex-coloring
$2$-edge-weighting. In \cite{CLWY}, Chang \textit{et al.} considered
vertex-coloring $2$-edge-weighting in bipartite graphs and proved
the following results.

\begin{lemma}{\rm{(Chang, Lu, Wu and Yu, \cite{CLWY})}} \label{yu}

Every connected nice bipartite graph admits a vertex-coloring
$2$-edge-weighting if one of following conditions holds:

    $(1)$ $|A|$ or $|B|$ is even;

    $(2)$ $\delta(G) =1$;

    $(3)$ $\lfloor d(u)/2 \rfloor +1 \neq d(v)$ for any edge $uv \in
    E(G)$.
\end{lemma}

\begin{theorem}\label{thm42} Let $G$ be a nice graph. If
$\delta(G)\geq 8\chi(G)$, then $G$ admits a vertex-coloring
$2$-edge-weighting.
\end{theorem}

\pf  Let $\{V_1,\ldots, V_{\chi(G)}\}$ be a partition of $V(G)$ into
independent sets. For each $v\in V_{i}$, choose $a_{v}^{-}$ such that
$\lfloor\frac{d(v)}{4}\rfloor\leq a_{v}^{-}\leq
\lfloor\frac{d(v)}{2}\rfloor$, $a_{v}^{-}+d_{G}(v)\equiv 2i$ (mod
2$\chi(G))$, and $a_{v}^{-}+2\chi(G)\geq
\lfloor\frac{d(v)}{2}\rfloor$. Such choice for $a_{v}^{-}$ exists as
$\delta(G)\geq 8\chi(G)$. Set $a^{+}_{v}=a_{v}^{-}+2\chi(G)$.

Furthermore, such a choices of $a_{v}^{-}$ and $a^{+}_{v}$ satisfy
the conditions of Lemma \ref{Berry}, i.e.,
\begin{align*}
\frac{1}{2}(d(v)-a_v^{-}-2\chi(G))-\chi(G)&=\frac{1}{2}(d(v)-a_v^{+})-\chi(G)\\
&\geq \frac{d(v)}{8}-\chi(G),
\end{align*}
thus there is a subgraph $H$ such that for all $v$,
$d_{H}(v)\in\{a_{v}^{-},a_{v}^{-}+1,a_{v}^{+},a_{v}^{+}+1\}$. Set
 $w(e)=2$ for $e\in E(H)$ and $w(e)=1$ for $e\in
 E(G)-E(H)$. If $v\in V_{i}$, we have
\begin{align*}
\sum_{v\sim e}w(e)=2d_{H}(v)+d_{G-H}(v)=d_{G}(v)+d_{H}(v)\in
\{2i,2i+1\} \ (\mbox{mod }\ 2\chi(G)).
\end{align*}
Thus adjacent vertices in different parts of $\{V_1,\ldots,
V_{\chi(G)}\}$ have different arities. As each $V_i$ is an
independent set, these weights form a vertex-coloring
2-edge-weighting of G. \qed

\vspace{2mm}

\begin{theorem}
Given a nice bipartite graph $G=(U,W)$, if there exists a vertex $v$
such that $d_{G}(v)\not \in \{d_{G}(x)\ |\ x\in N(v)\}$ and
$G-v-N(v)$ is connected, then $G$ admits a vertex-coloring
$2$-edge-weighting.
\end{theorem}

\pf If $|U|\cdot|W|$ is even, by Lemma \ref{yu}, the result follows.
So we may assume that both $|U|$ and $|W|$ are odd. Let $v\in U$
satisfy $d_{G}(v)\not \in \{d_{G}(x)\ |\ x\in N(v)\}$ and
$N(v)=\{w_{1},\ldots,w_{k}\}$. Since $G-v-N(v)$ is connected, by
Theorem \ref{group3}, $G-v-N(v)$ has a vertex-coloring
$2$-edge-weighting such that $c(x)$ is odd  for all $x\in U-v$ and
$c(y)$ is even for all $y\in W-N(v)$. Now we assign every edge of
$E[N(v),U]$ with weight $2$. Clearly $c(x)$ is odd for all $x\in
U-v$ and $c(y)$ is even for all $y\in W$. Moreover $c(v)\neq c(u)$
for all $u\in N(v)$ since $d(u)\neq d(v)$. Thus we obtain a
vertex-coloring $2$-edge-weighting of $G$. \qed

\begin{theorem}
Given a nice bipartite graph $G=(U,W)$, if there exists a vertex $v$
of degree $\delta(G)$ such that $d_{G}(v)\not \in \{d_{G}(x)\ |\
x\in N(v)\}$ and $G-v$ is connected, then $G$ admits a
vertex-coloring $2$-edge-weighting.
\end{theorem}

\pf If $|U|\cdot|W|$ is even, by Lemma \ref{yu}, the result follows.
So we may assume that both $|U|$ and $|W|$ are odd. Let $v\in U$
satisfy $d_{G}(v)=\delta(G)$ and $d_{G}(v)\not \in \{d_{G}(x)\ |\
x\in N(v)\}$ . Now we consider two cases.

{\it Case 1.} $\delta(G)$ is even.

In this case, $|(U-v)\cup N(v)|$ is even and $W-N(v)$ is odd. By
Theorem \ref{group3}, $G-v$ has a vertex-coloring $2$-edge-weighting
such that $c(x)$ is odd for all $x\in (U-v)\cup N(v)$ and  $c(y)$ is
even for all $y\in W-N(v)$. Since $d_{G}(v)\not \in \{d_{G}(x)\ |\
x\in N(v)\}$, assigning the edges incident to $v$ with weight $1$
induces a vertex-coloring $2$-edge-weighting of $G$.

{\it Case 2.} $\delta(G)$ is odd.

In this case, $|(U-v)\cup N(v)|$ is odd and $W-N(v)$ is even. By
Theorem \ref{group3}, $G-v$ has a vertex-coloring $2$-edge-weighting
such that $c(x)$ is even for all $x\in (U-v)\cup N(x)$ and  $c(y)$
is odd for all $y\in W-N(v)$. Since $d_{G}(v)\not \in \{d_{G}(x)\ |\
x\in N(v)\}$, assigning the edges incident to $v$ with weight $1$
induces a vertex-coloring $2$-edge-weighting of $G$.  \qed

\section{3-connected bipartite graphs}

An interesting corollary of Lemma \ref{yu} is that every {\it
$r$-regular} nice bipartite graph ($r\geq 3$) admits a
vertex-coloring $2$-edge-weighting. Note $C_6,C_{10},\ldots$ are
2-regular nice bipartite graphs which do not admit
vertex-coloring 2-edge-weightings. 

In the following, we continue the research in this direction and
prove that a vertex-coloring 2-edge-weighting exists for every
3-connected bipartite graph. The following lemma is an important
step in proving our main result.

\begin{lemma}\label{lu}
Let $G$ be a $3$-connected non-regular bipartite graph with
bipartition $(U,W)$. Let $u\in U$ with $d(u)=\delta(G)$ and $t\leq
\delta-1$. Denote $N^{\delta}(u)=\{v\ |\ d(v)=\delta,v\in
N_{G}(u)\}=\{u_{1},\ldots,u_{t}\}$. Then there exist
$e_{1},\ldots,e_{t}$, where $e_{i}$ is incident to vertex $u_i$ in
$G-u$ for $i=1,\ldots,t$, such that $G-u-\{e_{1},\ldots,e_{t}\}$ is
connected.
\end{lemma}

\pf Let $C_{1},\ldots,C_{s}$ be the components of
$G-u-N^{\delta}(u)$. We construct a bipartite multi-graph $H$ with
bipartition $(X,Y)$, where $X=\{u_{1}\ldots,u_{t}\}$,
$Y=\{c_{1},\ldots,c_{s}\}$ and
$|E_{H}(u_{i},c_{j})|=|E_{G}(u_{i},C_{j})|$ for $1\leq i\leq t$ and
$1\leq j\leq s$. Then $d_{H}(u_{i})=\delta-1$ for every $u_{i}\in
X$.

{\it Claim.} {\it  $H$ contains a connected spanning subgraph $T$
such that $d_{T}(v)\leq \delta-2$ for every $v\in X$.}

Suppose that the claim does not hold. Let $R$ be a connected induced
subgraph of $H$ satisfying

$i)$. $R$ contains a connected spanning subgraph  $M$ such that
$d_{M}(v)\leq \delta-2$ for every $v\in V(M)\cap X$;

$ii)$. $|V(R)|$ is maximum.

It is easy to see that $V(R)\neq\emptyset$ and $R\neq H$. Let
$R=(A,B)$, where $A\subseteq X$ and $B\subseteq Y$.  By maximality
of $R$,  we have $d_{R}(v)\geq \delta-2$ for every $v\in A$ and
$E_{H}(B,X-A)=\emptyset$. Let $L=\{v\ |\ d_{R}(v)=\delta-2,v\in
A\}$. We see $|L|\geq 2$ since $G$ is $3$-connected. Let $M^*$ be a
connected spanning subgraph of $R$ such that $d_{M^*}(v)=\delta-2$
for every $v\in A$. Note that for every connected spanning subgraph
$N^*$ of $M^*$, we have $d_{N^*}(w)=\delta-2$ for $w\in L$ by the
maximality of $R$. So every edge incident with $w$ in $M^*$, where
$w\in L$, is a cut-edge of $M^*$.  Let $|L|=l$ and
$|E(R)-E(M^*)|=m$. Then $l+m\leq t\leq \delta-1$. We have
\begin{align*}
\omega(M^*-L)=\omega(H-L-(E(R)-E(M^*)))-1\geq (\delta-3)l+1.
\end{align*}
So $\omega(H-L)\geq (\delta-3)l+2-m$, which implies
\begin{align*}
\omega(G-u-L)&\geq (\delta-3)l+1-m+1\\
&\geq (\delta-3)l+2-(\delta-1-l)\\
&=(\delta-2)l+3-\delta.
\end{align*}
Since $G$ is $3$-connected, then
\begin{align*}
3\omega(G-u-L)\leq (\delta-1)l+\delta-l.
\end{align*}
It follows that
\begin{align*}
\omega(G-u-L)\leq \lfloor\frac{(\delta-1)l+\delta-l}{3}\rfloor.
\end{align*}
However
\begin{align*}
&(\delta-2)l+3-\delta-\frac{(\delta-1)l+\delta-l}{3}=\frac{2\delta
l}{3}-\frac{4\delta}{3}-\frac{4l}{3}+3>0,
\end{align*}
a contradiction. So we complete the claim and thus obtain a
connected spanning subgraph $T$ of $H$.

Let $E'$ denote the set of corresponding edges of $E(T)$ in $G$.
Then we obtain a spanning subgraph $T^*=\bigcup_{i=1}^{s}C_{i}\cup
N^{\delta}(u)\cup E'$ of $G-u$ such that $d_{T^*}(v)\leq \delta-2$
for every  $v\in N^{\delta}(u)$. Thus the proof is complete. \qed

\vspace{3mm}

The following theorem is the main result of this section.

\begin{theorem} \label{3connected} Let  $G=(U,W)$ be a nice bipartite graph.
If $G$ is $3$-connected, then $G$ admits a vertex-coloring
$2$-edge-weighting.
\end{theorem}

\pf If $G$ is a regular graph,  the result follows from Lemma
\ref{yu} (3). In the following, let $G$ be a $3$-connected
non-regular bipartite graph with bipartition $(U,W)$. Let $u\in U$
with $d(u)=\delta(G)$ and $N^{\delta}(u)=\{v\ |\ d(v)=\delta,v\in
N_{G}(u)\}=\{u_{1},\ldots,u_{t}\}$, where $t\leq \delta-1$. Then by
Lemma \ref{lu}, there exist $e_{1},\ldots,e_{t}$, where $e_{i}$ is
incident to vertex $u_i$ in $G-u$ for $i=1,\ldots,t$, such that
$G-u-\{e_{1},\ldots,e_{t}\}$ is connected.

By Lemma \ref{yu}, we can assume that $|U||W|$ is odd. Now we
consider two cases.

{\it Case 1.} $\delta(G)$ is even.

Then $|N(u)\cup (U-u)|$ is even.  By Theorem \ref{group3},
$G-u-\{e_{1},\ldots,e_{t}\}$ has a vertex-coloring
$2$-edge-weighting such that $c(x)$ is odd for all $x\in N(u)\cup
(U-u)$ and $c(y)$ is even for all $y\in W-N(u)$.  We assign every
edge of  $\{e_{1},\ldots,e_{t}\}$ with weight $2$ and every edge of
$\{uu_{i}\ |\ i=1,\ldots,t\}$ with weight $1$. If $d(u_{i})=d(u)$
for $i=1,\ldots,t$, then $d_{G-u-\{e_{1},\ldots,e_{t}\}}(u_{i})$ is
even. Now $c(u)=d(u)$ and $c(u)< c(u_{i})$ for $i=1,\ldots, t$.
Moreover, $c(u_{i})$ is even for $i=1,\ldots, t$. Hence we obtain a
vertex-coloring $2$-edge-weighting of $G$.

{\it Case 2.} $\delta(G)$ is odd.

Then $|W-N(u)|$ is even. By Theorem \ref{group3},
$G-u-\{e_{1},\ldots, e_{t}\}$ has a vertex-coloring
$2$-edge-weighting such that $c(x)$ is even for all $x\in N(u)\cup
(U-u)$ and $c(y)$ is odd for all $y\in W-N(u)$. We again assign
every edge of $\{e_{1},\ldots,e_{t}\}$ with weight $2$ and every
edge of $\{uu_{i}\ |\ i=1,\ldots,t\}$ with weight $1$. Then
$c(u)=d(u)$ and $c(u)< c(u_{i})$ for $i=1,\ldots,t$. Moreover,
$c(u_{i})$ is odd for $i=1,\ldots,t$. Then we obtain a
vertex-coloring $2$-edge-weighting of $G$.

We complete the proof. \qed

Based on the proof of Theorem \ref{3connected}, we can easily obtain
the following corollary.

\begin{coro}
Let $G=(U,W)$ be a bipartite graph with $\delta(G)\geq 3$. If there
exists a vertex of degree $\delta(G)$ such that $G-u-N(u)$ is
connected, then $G$ admits a vertex-coloring $2$-edge-weighting.
\end{coro}

\section{Conclusions}

In this paper, we prove that every $3$-connected bipartite graph has
a vertex-coloring $2$-edge-weighting. There exists a family of
infinite bipartite graphs (e.g., the generalized $\theta$-graphs)
which is 2-connected and has a vertex-coloring $3$-edge-weighting
but not a vertex-coloring $2$-edge-weighting. It remains an open
problem to classify all $2$-connected bipartite graphs admitting a
vertex-coloring $2$-edge-weighting.

\vspace{2mm}

\noindent \title{\Large\bf Acknowledgments} \maketitle The authors
are indebted to Dr. Yinghua Duan for the valuable discussion.  The
authors also express gratitude to Nankai University, where this
research was carried out, for the support and the hospitality.

\end{document}